\def\E{\end{document}}
\documentclass[11pt]{article}
\usepackage{amsmath}
\usepackage{amsfonts}
\usepackage{mathrsfs}
\usepackage{amssymb}
\usepackage{amssymb,amsmath}

\allowdisplaybreaks
\begin{document}
\title{Nonsmoothable involutions on spin 4-manifolds\thanks{This work is supported  by NSFC (10771023 and 10931005)}}
\author{Changtao Xue and  Ximin Liu}
\date{}
\newtheorem{theorem}{Theorem}[section]
\newtheorem{definition}{Definition}[section]
\newtheorem{lemma}{Lemma}[section]
\newtheorem{proposition}{Proposition}[section]
\newtheorem{corollary}{Corollary}[section]
\newtheorem{remark}{Remark}
\renewcommand{\theequation}{\thesection.\arabic{equation}}
\catcode`@=11 \@addtoreset{equation}{section} \catcode`@=12
\maketitle{}
\begin{center}{\bf Abstract}\end{center}
\begin{center}
\begin{minipage}{125mm}
Let $X$ be a closed, simply-connected, smooth, spin 4-manifold
whose intersection form is isomorphic to $n(-E_8)\bigoplus mH$,
where $H$ is the hyperbolic form.
 In this paper, we prove that for $n$ such that $n\equiv 2
~{\rm mod} ~4$, there exists a locally linear pseudofree
$\mathbb{Z}_2$-action on $X$ which is nonsmoothable with respect
to any possible smooth structure on $X$.



\end{minipage}
\end{center}

\section{Introduction}

 A topological finite group $G$-action on
an $n$-dimensional manifold $X$ is called locally linear if for
any point $x\in X $, there exists  a $G_{x}$-invariant
neighborhood $V_{x}$ of $x$ such that $V_{x}$ is homeomorphic to
$\mathbb{R}^n$, and $G_{x}$ acts on $V_{x}$ in a linear orthogonal
way, where $G_{x}$ is the isotropy group of $x$.

It is well-known that every smooth action is locally linear. On
the other hand, a locally linear action is not necessarily smooth.
In [8], S. Kwasik and T. Lawson proved the existence of
nonsmoothable actions on certain contractible 4-manifolds mainly
by gauge theory, and in some cases of involutions, Rohlin's
$\mu$-invariant is used. In recent years, many nonsmoothable group
actions on 4-manifolds are constructed by many authors [2, 3, 7,
9, 10, 11] including the second author. Almost all of them use
gauge theory to prove the nonsmoothability except the case of
nonsmoothable involution on $K3$ in [11] by Nakamura, where he
just uses the $G$-spin theorem to prove the nonsmoothability.

In this paper we restrict our attention to involutions on a class
of  spin 4-manifolds and prove the following theorem.

{\bf Theorem 1.1.} {\it Let $X$ be a closed, simply-connected,
smooth, spin 4-manifold
 whose intersection form is isomorphic to $n(-E_8)\bigoplus mH$, where $H$
is the hyperbolic form. If $n\equiv 2~ {\rm mod}~4$, then there
exists a locally linear pseudofree $\mathbb{Z}_2$-action on $X$
which is nonsmoothable with respect to any possible smooth
structure on $X$.}

A group action is said to be pseudofree if each nontrivial group
element has a discrete fixed point set.

The proof of  Theorem 1.1 is divided into two steps. In the first
step, we give a constraint on smooth involutions. In the second
step, we construct a locally linear action which would violate the
constraint.

To obtain a constraint on smooth involutions, we do not use gauge
theory, but only use Rohlin's theorem. On the other hand, to
construct a locally linear action, we use the realization theorem
due to  Edmonds and Ewing [4]. In fact, the method of the proof of
Theorem 1.1 is essentially same with Nakamura's method for the
construction of a nonsmoothable involution on $K3$ in [11]. While
Nakamura uses the $G$-spin theorem instead of Rohlin's theorem,
these two methods can be seen essentially same by recalling the
proof of Rohlin's theorem by the index theorem.
 By Freedman's theorem [5],
the homeomorphism type of a simply-connected 4-manifold is uniquely
determined by its intersection form if the intersection form is
even. Thus, our result can be applied to many spin 4-manifolds
including certain elliptic surfaces. Let $X=E(n)$ be the relatively
minimal elliptic surface with rational base. Then the intersection
form of X is isomorphic to $n(-E_{8})\oplus (2n-1)H$. So from
Theorem 1.1, we can get the following corollary immediately.

{\bf Corollary 1.2.}  {\it Let $X=E(n)$ be the relatively minimal
elliptic surface with rational base. If $n\equiv 2~ {\rm mod}~4$,
 then there exists a locally linear pseudofree
$\mathbb{Z}_2$-action on $X$ which is nonsmoothable with respect
to any possible smooth structure on $X$.}

In [7],  K. Kiyono  proved that if $X$ is a closed,
simply-connected, spin topological 4-manifold not homeomorphic to
either $S^{4}$ or $S^2 \times S^2$, then for any sufficiently large
prime number p, there exists a homologically trivial, pseudofree,
locally linear action of $\mathbb{Z}_{p}$ on $X$ which is
nonsmoothable. So it would be interesting to compare our result with
Kiyono's result, since we have constructed nonsmoothable
$\mathbb{Z}_p$-actions of the minimal order $p=2$ on a large class
of spin 4-manifolds.

The paper is organized as follows. In section 2, we provide some
preliminaries. In particular, a constraint on smooth involutions is
given. In section 3, we prove our main theorem (Theorem 1.1).

{\bf Acknowledgments:} The authors would like to thank referee for
many valuable suggestions, especially for pointing out that the
main theorem of this paper can be extended to a larger class of
spin 4-manifolds. It is also pleasure to thank N. Nakamura for his
help, in particular, for explaining how to prove Proposition 2.6.

\section{Preliminaries}

In this section, a constraint on smooth involutions and some
method of constructing locally linear $\mathbb{Z}_{2}$-actions are
given. Note that this section largely depends on Nakamura's paper
[11], so we refer the reader to the Nakamura's excellent
exposition [11] for more details.

{\bf 2.1 Construction of locally linear $\mathbb{Z}_2$-actions}

 To construct locally linear $\mathbb{Z}_2$-actions, we use the following special
case of the realization theorem by Edmonds and Ewing [4].

 {\bf Theorem 2.1} ([4, 11]). {\it Suppose that we are given a
$\mathbb{Z}_2$-invariant bilinear unimodular even form $\Psi :
V\times V\rightarrow \mathbb{Z}$} {\it which satisfies the following:\\
(1) As a $\mathbb{Z}[\mathbb{Z}_2]$-module, $V\cong T\oplus F$,
where $T$ is a trivial ${\mathbb{Z}}[\mathbb{Z}_2]$-module with
$rank_{\mathbb{Z}}
T=n$, and $F$ is a free ${\mathbb{Z}}[\mathbb{Z}_2]$-module.\\
(2) For any $v\in V$, $\Psi (v,gv)\equiv 0 ~{\rm mod}~ 2$, where g
is the
generator of $\mathbb{Z}_2$.\\
(3) The $G$-signature formula is satisfied, i.e.,
$\sigma(g,(V,\Psi))=0$.\\
Then, there exists a locally linear $\mathbb{Z}_2$-action on a
simply-connected 4-manifold $X$ such that its intersection form is
$\Psi$, and the number of fixed points is $n+2$.}\\

Since the form $\Psi$ is assumed even, the homeomorphism type of $X$
is unique by Freedman's theorem [5].

 For our application, we also need their  equivariant handle construction.

Let $B_{0}$ be a unit ball in $\mathbb{C}^2$, on which
$\mathbb{Z}_2$ acts by multiplication of $\pm 1$. Take a
$\mathbb{Z}_2$-invariant knot $K$ in $S_{0}=\partial B_{0}$. Then
 a framing  $r$ of  $K$  can  be represented by an equivariant embedding
$f_r: \, S^{1}\times D^{2}\rightarrow S_{0}$ for some
$\mathbb{Z}_2$ action on $S^{1}\times D^{2}$  given by
$g(z,w)=(-z,(-1)^{r-1}w)$. So for a given $K$ and a framing $r$, a
4-manifold $W=B_{0}\cup _{f_{r}}D^{2}\times D^{2}$ with a
$\mathbb{Z}_2$-action can be constructed.

Let $H_1, \cdots, H_n$ be copies of $D^2\times D^2$ on which $Z_2$
acts by $g(z,w)=(-z,-w)$. If a knot $K$ with an even framing $r$ is
given,  then we can attach $H_{i}$ to $B_0$ equivariantly via
$f_{r}$.

We will choose a $\mathbb{Z}_2$-invariant $n$-component framed link
$L$ in $\partial B_{0}$ to represent $\Psi$ as follows. Under the
assumption of Theorem 2.1, we may assume $\Psi |_{T}$ is represented
by a matrix $(a_{ij})$ such that $a_{ii}$ is even and $a_{ij}$ is
odd whenever $i\neq j$. Actually, we can take a $n$-component link
$L_{T}$ in $S_{0}=\partial B_0$ representing the matrix $(a_{ij})$
such that each component of $L_T$ is $\mathbb{Z}_2$-invariant. And
it is not difficult to realize the other part of $\Psi$ by a link,
and therefore we obtain a framed link $L$ in $S_0$ which realizes
the given $\mathbb{Z}_2$-invariant form $\Psi$. Attach
$H_{1},\ldots, H_{n}$ and free 2-handles to $B_{0}$ along $L$
equivariantly. Thus we obtain a 4-manifold $X_{0}$ on which
$\mathbb{Z}_2$ acts smoothly,
$$X_{0}=B_{0}\cup H_{1}\cup \cdots \cup H_{n}\cup({\rm free~ handles}).$$
The boundary of $X_{0}$ is an integral homology 3-sphere $\Sigma$
with a free $\mathbb{Z}_2$-action. Edmonds and Ewing proved that
there exists a contractible 4-manifold $Z$ with a locally linear
$\mathbb{Z}_{2}$-action such that its boundary is $\Sigma$ with
the given free $\mathbb{Z}_{2}$-action, and it has exact one fixed
point. Then we obtain the required manifold $X=X_{0}\cup Z$ with
the required action.

 Note that each of the components $B_{0}$, $H_{1}$, $\cdots$, $H_{n}$, $W$ has one
fixed point, denoted by $P, Q_{1}, \cdots, Q_{n}, P^{'}.$ The action
constructed above is smooth on $X_{0}$, and is smooth on $X$ except
near the final fixed point.

{\bf 2.2 Atiyah-Bott's criterion for $\varepsilon(P)$}

For a smooth even-type (pseudofree equivalently)
$\mathbb{Z}_2$-action on a simply-connected smooth spin 4-manifold
$X$ which preserves the unique spin structure and also the ${\rm
Spin}^{c}$-structure $c_0$ which is determined by the spin
structure. The sign assignment $\varepsilon$ is introduced by Atiyah
and Bott in [1]. They defined the Spin-number to be the Lefschetz
number of the corresponding Dirac complex. The Spin-number has the
form ${\rm ind}_{g}D=\sum _{P\in {X^{\mathbb{Z}_2}}}\nu(P)$, where
$g$ is the generator of the $\mathbb{Z}_2$-action. The sign of each
summand $\nu(P)$ depends on P and $g$, that is the assignment
$\varepsilon(P)$. In other worlds, the sign assignment determined by
the lift of the $\mathbb{Z}_2$-action to $c_0$ is $\varepsilon :
X^{\mathbb{Z}_2}\rightarrow {\pm1}$. By the $G$-spin theorem [1], we
have $${\rm ind}_{g}D=k_{+}-k_{-}=\frac{1}{4}\sum _{P\in
{X^{\mathbb{Z}_2}}}\varepsilon(P),$$ $${\rm
ind}D=k_{+}+k_{-}=-\frac{1}{8}\sigma(X),$$ where $k_{+}$ and $k_{-}$
are coefficients of the $\mathbb{Z}_2$-index of the Dirac operator.
Note that $k_{+}$ and $k_{-}$ are even because of the quartanionic
structure of Dirac index. Then, we can see that the sum $\sum_{P\in
X^{\mathbb{Z}_2}} \varepsilon(P)$ is a multiple of 8 by solving the
above equations. Suppose $g\in \mathbb{Z}_2$ is nontrivial element,
then we can lift the smooth pseudofree involution $g: X\rightarrow
X$ to the frame bundle $F$ as $g_{*}: F\rightarrow F$ if a
$g$-invariant metric is fixed. A spin structure on $X$ is given by a
double cover $\varphi: \hat{F}\rightarrow F$, where $\hat{F}$ is a
Spin(4)-bundle. Suppose that $g$ lifts to $\hat{F}$. The values
$\varepsilon(P)$ and $\varepsilon(Q)$ of distinct fixed points $P$
and $Q$ can be compared by Atiyah-Bott's criterion as following
proposition. We take a path $s$ in $F$ starting form a point $y\in
F_{P}$ and ending at $y^{'}\in F_{Q}$. Then the path $-g_{*}s$ has
the same starting point and the end point as $s$, where $"-"$ means
the multiplication by $-1$ on each fiber. Thus by connecting $s$ and
$-g_{*}s$, we obtain a circle $C$ in $F$.

{\bf Proposition 2.2} ([1, 11]). {\it The preimage
$\varphi^{-1}(C)$ has two components if and only if
$\varepsilon(P)=\varepsilon(Q)$. In other words, the preimage
$\varphi^{-1}(C)$ is connected if and only if
$\varepsilon(P)=-\varepsilon(Q)$.}

Recall that each component of $B_0$ and $H_{i}$ of $X_0$ constructed
in subsection 2.1 has a fixed point, denoted by $P$ and $Q_i$,
suppose that $H_i$ is attached to $B_0$ equivariantly along a knot
$K$ with a framing $r$. Then we have the following proposition.

{\bf Proposition 2.3} ([11]). {\it Suppose $K$ is a trivial knot
in $\partial B_0$ which bounds a $\mathbb{Z}_2$-invariant embedded
disk $D_0$ in $B_0$ containing $P$. If $r\equiv 2 ~{\rm mod}~4$,
then $\varepsilon(P)=\varepsilon(Q)$. If $r\equiv 0 ~{\rm mod}~4$,
then $\varepsilon(P)=-\varepsilon(Q)$.}

Let X be an oriented topological manifold and let
$G$=$\mathbb{Z}_2$. The sign assignment can be defined for locally
linear actions by using Atiyah-Bott's criterion itself on
topological spin structure $\varphi: \hat{F}\rightarrow F$, and it
depends only on equivalence classes of orientation-preserving
locally linear $G$-actions on $X$. See [11] for more details. On
the fiber of $F$ over each fixed point $P$, there is a point
$y_{P}$ which is mapped to $-y_{P}$ by the $\mathbb{Z}_2$-action.
For distinct fixed points $P$ and $Q$, by taking a path $s$
connecting such a $y_{P}$ with such $y_{Q}$, we can define the
sign assignment.

{\bf Definition 2.4} ([11]). {\it For each pair ($P,~Q$) of fixed
points, let $s$ be a path in $F$ as above, and $C$ the circle
formed by $s$ and $-g_{*}s$. Define $\varepsilon^{'}(P,Q)$ by
$\varepsilon^{'}(P,Q)=1, ~if ~\varphi^{-1}(C)~ has~ 2
~components$;
$\varepsilon^{'}(P,Q)=-1, ~if ~\varphi^{-1}(C)~ is~ connected.$}\\

Note that this definition does not depend on smooth structures,
and is well-defined if $X$ is simply-connected. Furthermore, if
the action is realized by a smooth action, then
\begin{align}
\varepsilon^{'}(P,Q)=\varepsilon(P)\varepsilon(Q).
\end{align}

{\bf 2.3 Constraint on smooth involutions}

 We give a constraint on
smooth involutions by considering Rohlin's theorem. Rohlin's
theorem gives a criteria for the question that which
simply-connected topological manifolds carry smooth structures.

{\bf Theorem 2.5} ([5]). {\it If X is a smooth, closed, spin
4-manifold, then the signature $\sigma(X)$ of $X$ satisfies
$\sigma(X)\equiv 0~{\rm mod}~16$. }

Now let $X$ be a smooth, closed, oriented, simply-connected  spin
4-manifold, and suppose that $\mathbb{Z}_{2}$ acts on $X$ smoothly
and pseudofreely in an orientation-preserving way. By [1], the
fixed-point set $X^{\mathbb{Z}_2}$ is discrete if and only if the
$\mathbb{Z}_{2}$-action lifts to the spin structure. Therefore $X
/ \mathbb{Z}_2$ is a spin $V$-manifold. The quotient singularities
are cones of $\mathbb{R}P^3$. It is well known that the spin
structures on $\mathbb{R}P^3$ can be divided into two equivalent
classes $s_{\pm}$ which are characterized as follows: let
$\widetilde{s}_{\pm}$ be the unique spin structure on the disk
bundle $D_{\pm}$ over $S^{2}$ of degree $\pm 2$, then
$s_{\pm}=\widetilde{s}_{\pm}|_{\partial D_{\pm}}$. Define the spin
type of a fixed point by the spin structure on $\mathbb{R}P^3$
induced from $X / \mathbb{Z}_2$. The number of fixed points
corresponding to $s_{\pm}$ is denoted by $n_{\pm}$. We know that
the quotient space $X/\mathbb{Z}_2$ is not a smooth 4-manifold. In
order to use Rohlin's theorem, we need to make $X/\mathbb{Z}_2$ to
be smooth. Since $\mathbb{R}P^3$ has two equivalent classes of
spin structures, we need to make the spin structures are
compatible. Remove cones of $\mathbb{R}P^3$ from $X /
\mathbb{Z}_2$, and glue disk bundles $D_{+}$ and $D_{-}$ so that
spin structures are compatible. We get a smooth spin 4-manifold.
Then, applying Rohlin's theorem, we have the following formula:

$$\sigma(X / \mathbb{Z}_2)\equiv n_{+}-n_{-}{\rm mod} ~16.$$
 Together with the $G$-signature theorem [1], we have
 $$\frac{1}{2}\sigma(X)\equiv n_{+}-n_{-} {\rm mod} ~16.$$
 So for a smooth
$\mathbb{Z}_{2}$-action on $X$ as above, we have
\begin{equation} \label{eq:2}
\left\{ \begin{aligned}
         \sharp  X^{\mathbb{Z}_2} &= n_{+}+n_{-}, \\
              \frac{1}{2}\sigma(X)&\equiv n_{+}-n_{-}{\rm mod}
              ~16.
                          \end{aligned} \right.
                          \end{equation}

There is a relation between the spin types and the sign assignments
of two fixed points.

 {\bf Proposition 2.6.} {\it Let $P$ and $Q$ be
distinct fixed points. Then, $\varepsilon(P)$=$\varepsilon(Q)$ iff
$P$ and $Q$ have the same spin type;
$\varepsilon(P)$=$-\varepsilon(Q)$ iff $P$ and
$Q$ have the different spin types.}

{\bf Proof.} The quantity $\varepsilon$ was introduced in order to
see the lifting way of the $\mathbb{Z}_2$-action around each fixed
point. Actually, there are two lifting ways. We will explain it in
the following.

Let us consider a local model around a fixed point. Let $D$ be a
unit 4-ball in $\mathbb{R}^4$, and let $\mathbb{Z}_2$ acts on $D$
by $x\mapsto -x$. The frame bundle of $D$ can be assumed as a
trivial bundle $F=X \times {\rm SO(4)}$. The induced
$\mathbb{Z}_2$-action on $F$ is given by $(x,v)\mapsto (-x,-v)$
for $(x,v)\in X\times {\rm SO(4)}$. Recall that ${\rm
Spin(4)=Sp(1)\times Sp(1)}$, and there is a canonical homomorphism
$\phi_0:{\rm Spin(4)}\mapsto {\rm SO(4)}$. Let $\hat{F}=X\times
{\rm Spin(4)}$. Then the unique spin structure on $F$ is given by
$\phi_0: \hat{F} \mapsto  F$ defined by $\phi=id_{X}\times
\phi_0$. We can see that there are two ways of lifting the
$\mathbb{Z}_2$-action to $\hat{F}$. They are given by involutions
$\rho_1$ and $\rho_2$ as follows:$$\rho_1:(x, q_{+},q_{-})\mapsto
(-x, -q_{+},q_{-}),$$ $$\rho_2:(x, q_{+},q_{-})\mapsto (-x,
q_{+},-q_{-}).$$
 On $S^3=\partial D$, the $\mathbb{Z}_2$-action is
the anti-podal map. Then $\hat{F}_1\triangleq \hat{F}/\rho_1$ and
$\hat{F}_2\triangleq \hat{F}/\rho_2$ give two spin structures on
$S^3/\mathbb{Z}_2=\mathbb{R}P^3$. Since
$H^1(\mathbb{R}P^3;\mathbb{Z}_2)=\mathbb{Z}_2$, $\mathbb{R}P^3$
has two distinct isomorphism classes of spin structures. By the
construction, the difference of $ \hat{F}_1$ and $\hat{F}_2$ is
given by a nontrivial $\mathbb{Z}_2$-bundle on $\mathbb{R}P^3$,
which gives a nonzero class in
$H^1(\mathbb{R}P^3;\mathbb{Z}_2)=\mathbb{Z}_2$. So, the
isomorphism classes of $ \hat{F}_1$ and $\hat{F}_2$ are different.

Recall that the singularities of our $\mathbb{Z}_2$-actions are
cones of $\mathbb{R}P^3$. Now, consider the definitions of the
Spin types and the sign assignments of fixed points together with
the analysis above. The result follows.

For the fixed points obtained in the subsection 2.1, we can compare
their spin
  types by the following proposition. The proposition is obtained by
  considering Proposition 2.3 and Proposition 2.6 simultaneously.

{\bf Proposition 2.7.} {\it Suppose $K_{i}$ bounds a
$\mathbb{Z}_2$-invariant embedded disk in $B_0$, and $r_i$ is the
framing of $K_i$. Then, $r_{i}\equiv 2 ~{\rm mod}~4$ if and only if
$P$ and $Q_i$ have the same spin types; $r_{i}\equiv0~ {\rm mod}~4$
if and only if $P$ and $Q_i$ have the different spin types. }\\

We obtained a constraint on smooth involutions as equations (2.2).
By the following proposition, we can see that the nonsmoothablity
of the $\mathbb{Z}_2$-action in Theorem 1.1 does not depend on
smooth structures.

 {\bf Proposition 2.8.} {\it The numbers $n_{+}$ and $n_{-}$ do not depend on smooth structures, i.e. they
  are invariants of locally linear pseudofree involutions on topological spin 4-manifolds.}

{\bf Proof.}  We know that Definition 2.4 extends the sign
assignment $\varepsilon$
 to the case of locally linear $\mathbb{Z}_2$-actions on a
simply-connected topological 4-manifold $X$ with isolated fixed
points. And this definition does not depend on smooth structures.
Now, by the relation (2.1) and Proposition 2.6, we have the
result.

\section{Proof of Theorem 1.1}

~~~~~~Let X be a closed, simply-connected, smooth, spin 4-manifold
which has the intersection form isomorphic to $n(-E_8)\bigoplus
mH$, where $n$ being a positive integer satisfying $n~\equiv ~ 2~
{\rm mod}~ 4$. Therefore $m\geq n+1$ by Furuta's inequality [6].
Suppose an orientation-preserving smooth pseudofree
$\mathbb{Z}_2$-action on X is given.

Since $\sigma(X)=-8n$, so by (2.2), $-4n\equiv~ n_{+}-n_{-} ~{\rm
mod}~ 16$. Therefore, if $n_{+}=n_{-}$, then $n\equiv~ 0 ~{\rm
mod}~ 4$.

 Now we can
construct a locally linear $\mathbb{Z}_{2}$-action on $X$ following
the method in section 2. By Theorem 2.1, if we fix an appropriate
$\mathbb{Z}_{2}$-action on the intersection form, then we have a
locally linear $\mathbb{Z}_{2}$-action on X.

Define $\mathbb{Z}_{2}$-action on $\Psi$=$n(-E_{8})\bigoplus mH$ as
follows. Let us select a positive integer $r$ less than $m$ and such
that $r\equiv m~{\rm mod}~2$. Let $\mathbb{Z}_{2}$ act on a
$\frac{m-r}{2}H\bigoplus \frac{m-r}{2}H$ summand by permutation of
two  $\frac{m-r}{2}H$'s. Similarly, let $\mathbb{Z}_{2}$ act on a
$n(-E_8)=\frac{n}{2}(-E_{8})\bigoplus\frac{n}{2}(-E_{8})$ summand by
permutation of two $\frac{n}{2}(-E_{8})$'s, and on the rest $rH$
trivially. The trivial part is denoted by $T$.

Now  consider the matrix
$$A=\left(
\begin{array}{cccccc}
0& \cdots &1 & 1&\cdots &1\\
\vdots & \ddots& \vdots&  \vdots&\cdots &\vdots\\
 1& \cdots& 0 &  1&\cdots&1\\
1& \cdots& 1& 2&  \cdots& 1\\
\vdots&  \cdots& \vdots& \vdots&  \ddots&\vdots\\
1&  \cdots& 1&  1& \cdots& 2
\end{array}
\right),$$ which is  a $2r\times2r$ matrix such that the first $r$
diagonal entries are 0, and the rest $r$ diagonal entries are 2, and
all of off-diagonal entries are 1.

We claim that the symmetric form represented by $A$ is isomorphic to
$rH$. It is
 clear that the form is even and indefinite. By diagonalizing
 the matrix $A$, we can see that the numbers of $+1$ and $-1$ on the
 diagonal are equal. The determinant of $A$ is $(-1)^r$. So the claim
 follows.
Hence, we may assume $\Psi|_{T}$ is represented by the matrix $A$.
Furthermore, the matrix $A$ can be realized by a link whose every
component bounds a $\mathbb{Z}_2$-invariant embedded disk as
follows. Let $p:S^{3}\rightarrow S^{2}$ be the Hopf fibration. Then
the inverse image of distinct $2r$ points in $S^{2}$ by $p$ forms a
required link. As in section 2,  by equivariant handle construction,
we can construct a smooth action on a manifold $X_{0}$ with
boundary,  and by Theorem 2.1, this action can be extended to the
whole $X$ as a locally linear action.

By now we can construct a pseudofree, locally linear
$\mathbb{Z}_{2}$-action on $X$ which satisfies
 $$\sharp X^{\mathbb{Z}_{2}}=2r+2.$$
 Denote the fixed points by $P,~Q_1, \cdots,~Q_{2r},~P^{'}$ as in
 section
 2. Recall that the $\mathbb{Z}_{2}$-action we have constructed is smooth
 except near $P^{'}$.
 By Proposition 2.7, a half number of fixed points in $\bigcup_{1\leq i \leq 2r} Q_{i}$ have
 the same spin type with $P$, and the rest fixed points in $\bigcup_{1\leq i \leq 2r} Q_{i}$
 have different spin type with $P$. If the $\mathbb{Z}_{2}$-action can be smooth on $X$, by the fact that
 $\sum_{P\in X^{\mathbb{Z}_2}}\varepsilon(P)$ is a multiple of 8 and Proposition 2.3, we
 have $\sum_{P\in X^{\mathbb{Z}_2}}\varepsilon(P)=0$. Thus, $P$ and
 $P^{'}$ have the different spin types. By Proposition 2.6, we can see that
 $$(n_{+},n_{-})=(r+1,r+1).$$
 But for $X$ whose intersection form is isomorphic to $n(-E_8)\bigoplus
mH$ with $n\equiv 2~ {\rm mod} ~4$, we can not have
 $n_{+}=n_{-}$  by the relation (2.2). Thus Theorem 1.1 is proved.

{\bf Remark 3.1} Note that for the positive integer $n$ such that
$n\equiv 0 ~{\rm mod}
 ~4$, by the consideration in subsection 2.3, we have $$n_{+}-n_{-}\equiv 0 ~{\rm mod} ~16.$$
 The pairs $(n_{+},n_{-})$ with $n_{+}=n_{-}$ satisfy it, so we can not conclude the nonsmoothablity
 as the case of $n\equiv 2 ~{\rm mod}~ 4$. It is an interesting
 question that whether there exists a smooth involution on a closed, simply-connected, smooth, spin 4-manifold
which has the intersection form isomorphic to $n(-E_8)\bigoplus mH$
with $n\equiv 0 ~{\rm mod}~ 4$ and $n_{+}=n_{-}$. But we cannot
construct such an example now.

{\small
 \vspace{3mm}

\begin{flushleft}
School of Mathematical Sciences, Dalian University of Technology, Dalian 116024, China\\
E-mail: xuechangtao123@163.com\\
\vspace{2mm} Department of  Mathematics, South China University
of Technology, Guangzhou 510641, China\\
E-mail: ximinliu@scut.edu.cn
 \end{flushleft}

\end{document}